\documentclass[11pt]{amsart}
\usepackage{amssymb,amsmath,amsthm}
\newcommand{\leg}[2]{\genfrac{(}{)}{}{}{#1}{#2}}

\newtheorem{theorem}{Theorem}
\newtheorem{lemma}[theorem]{Lemma}
\newtheorem{corollary}[theorem]{Corollary}

\newtheorem{prop}[theorem]{Proposition}

\theoremstyle{remark}

\newtheorem*{remark}{Remark}

\numberwithin{theorem}{section} \numberwithin{equation}{section}

\newcommand{\calD}{\mathcal{D}}

\newcommand{\R}{\mathbb{R}}
\newcommand{\C}{\mathbb{C}}

\newcommand{\Q}{\mathbb{Q}}

\newcommand{\Z}{\mathbb{Z}}
\newcommand{\N}{\mathbb{N}}

\newcommand{\lcm}{{\text {\rm lcm}}}

\def\H{\mathbb{H}}

\newcommand{\beqs}{\begin{equation*}}
\newcommand{\eeqs}{\end{equation*}}
\newcommand{\beq}{\begin{equation}}
\newcommand{\eeq}{\end{equation}}
\newcommand\mylabel[1]{\label{#1}}
\newcommand\eqn[1]{(\ref{eq:#1})}

\newcommand\propo[1]{\ref{propo:#1}}
\newcommand\Etwid{\overset {\text{\lower 3pt\hbox{$\sim$}}}E}
\newcommand\Ftwid{\overset {\text{\lower 3pt\hbox{$\sim$}}}F}
\newcommand\Qtwid{\overset {\text{\lower 3pt\hbox{$\sim$}}}Q}
\newcommand\thm[1]{\ref{thm:#1}}
\newcommand\corol[1]{\ref{cor:#1}}

\begin{document}
\title[Quasiweak Maass forms]{Partition statistics and quasiweak
Maass forms}

\author{Kathrin Bringmann}
\address{School of Mathematics\\University of
Minnesota\\ Minneapolis, MN 55455 \\U.S.A.}
\email{bringman@math.umn.edu}
\author{Frank Garvan}
\address{Department of Mathematics, University of Florida, Gainesville,
FL 32611-8105, U.S.A.}
\email{frank@math.ufl.edu}
\author{Karl Mahlburg}
\address{Massachusetts Institute of Technology
Department of Mathematics
MA 02139-4307\\ U.S.A.}
\email{mahlburg@math.mit.edu}
\thanks{The second author was supported in part by NSA Grant H98230-07-1-0011.  The third author 
was partially supported by a Clay Liftoff Fellowship.}

\date{February 19, 2008}

\begin{abstract}

Andrews recently introduced $k$-marked Durfee symbols, which are a generalization of partitions
that are connected to moments of Dyson's rank statistic.  He used these connections to find
identities relating their generating functions as well as to prove Ramanujan-type 
congruences for these objects and find relations between.

In this paper we show that the hypergeometric generating functions for these objects are natural
examples of {\it quasimock theta functions}, which are defined as the holomorphic parts of weak Maass forms and
their derivatives.  In particular, these generating functions may be viewed as analogs of Ramanujan's mock theta functions with arbitrarily high weight.  We use the automorphic properties to prove the existence of 
infinitely many congruences for the Durfee symbols.  Furthermore, we show that as $k$ varies, the modularity of the $k$-marked Durfee symbols is precisely dictated by the case $k=2.$  Finally, we use this relation in order to prove the existence of general congruences for rank moments in terms of level one modular forms of bounded weight.

\end{abstract}
\maketitle

\section{Introduction and Statement of results}

Modular and automorphic forms play an important role in many different
areas, including mathematical physics, representation theory, the theory of elliptic curves, quadratic forms,
and partitions, just to mention a few.  Many important generating
functions are modular forms, and the modular transformation properties can often be
used to prove arithmetic properties for the underlying combinatorial objects \cite{Ono}.  
There are many other examples of modular forms found in the realm of
hypergeometric $q$-series, such as the infinite products in the Rogers-Ramanujan identities.  
Until recently, however, there were few known examples of more general 
automorphic forms arising from similar ``arithmetic" generating functions.  Work of the first author and Ono
\cite{BO1, BO2}
(see also  \cite{B1,BFO,BL,BOR})  constructed several infinite
families of weak Maass forms of weights $1/2$ and $3/2$, whose
holomorphic parts were based on Ramanujan's mock theta functions and
also more general hypergeometric functions.

With the benefit of retrospect, we may view the error terms in the
mock theta transformations as suggesting that the appropriate
functions to consider were not modular forms, but the more general
weak Maass forms.  In this paper we study higher weight analogs that
are based on Andrews' work on Durfee symbols \cite{An2}.  The main
result of this paper is that the associated
generating functions  can be written in terms of automorphic
functions of higher weights, namely the derivatives of weak Maass
forms.  Although these derivatives are difficult to understand on
their own, we show that there are cancellations among certain linear combinations of
derivatives of different orders, and are able to successfully describe the analytic
behavior of the moment functions.  Furthermore, this allows us to
better understand the arithmetic of the coefficients; as a sample application we 
prove the existence of congruences.

In order to define the new objects at hand, first recall the
generating function for the partition function,
\begin{eqnarray}\mylabel{eq:partgen}
P(q)=P(z):= \sum_{n=0}^{\infty} p(n)\, q^{n}
= q^{\frac{1}{24}}\eta(z)^{-1},
\end{eqnarray}
where
$\eta(z):=
q^{ \frac{1}{24}}\, \prod_{n=1}^{\infty}(1-q^n)$
is Dedekind's $\eta$-function, a weight $\frac12$ modular form,  and
$q:=e^{2 \pi i z}$.
Of the the many consequences of the modularity properties of $P(q)$, some of the most striking are
the three congruences due to Ramanujan,
namely
\begin{align*}
p(5n+4) & \equiv 0 \pmod 5,\\
p(7n+5) & \equiv 0 \pmod 7, \\
p(11n+6) & \equiv 0 \pmod {11}.
\end{align*}

To explain the congruences with modulus $5$ and $7$, Dyson \cite{Dy}
introduced the \textit{rank} of a partition, which is defined to be
its largest part minus the number of its parts.
Dyson conjectured that the partitions of $5n+4$ (resp. $7n+5$) form
$5$ (resp. $7$) groups of equal size when sorted by their ranks
modulo $5$ (resp. $7$).
This conjecture was proven by Atkin and Swinnerton-Dyer \cite{AS}.
If $N(m,n)$ denotes the number of partitions of $n$  with rank $m$,
then
we have the generating  function
\begin{eqnarray*}
R(w;q):= 1 + \sum_{m \in \Z}
\sum_{n=1}^{\infty}
N(m,n) w^m q^n
=1 + \sum_{n=1}^\infty \frac{q^{n^2}}{(wq;q)_n (w^{-1}q;q)_n}
= \frac{(1-w)}{(q;q)_{\infty}} \sum_{n \in \Z}
\frac{(-1)^n q^{\frac{n}{2}(3n+1) }}{1-wq^n},
\end{eqnarray*}
where $(a;q)_{n}:= \prod_{j=0}^{n-1}(1-aq^j)$ and $(a;q)_{\infty}:=
\lim_{n \to \infty} (a;q)_n $.
In particular
\begin{eqnarray*}
R(1;q)&=&\mathcal{P}(q), \\
R(-1;q)&=&f(q):=1+\sum_{n=1}^{\infty}\frac{q^{n^2}}{(-q;q)_n^2 }.
\end{eqnarray*}
The function $f(q)$ is one of the \textit{mock theta functions}
defined by Ramanujan in his last letter to Hardy.  The first author
and Ono shed light on their mysteries by showing that if $w$ is a root of
unity, then the rank generating functions $R(w;q)$ (and in particular
$f(q)$)
are the ``holomorphic parts'' of  weak Maass forms \cite{BO2} (we say
more on these results and recall the notion of a weak Maass form in
Section \ref{MaassSection}).
The theory of weak Maass forms proved to be very useful for understanding the arithmetic of
the coefficients, leading to many notable results.  These include, for example, an exact formula
for the coefficients of $f(q)$ \cite{BO1}, asymptotics for
$N(m,n)$ \cite{B2},
identities for rank differences \cite{BOR}, and  congruences for
certain partition statistics   \cite{BO2}.

Here we consider infinite families of weak Maass forms of arbitrarily
high half-integer weight that also arise from combinatorial hypergeometric functions.
To state those results recall that Andrews introduced in
\cite{An2}  the    \textit{symmetrized $k$-th rank moment function}
\begin{eqnarray}
\eta_k(n):= \sum_{m = - \infty}^{\infty}
\left(
\begin{matrix}
m + \left[\frac{k-1}{2} \right] \\
k
\end{matrix}
\right)
N(m,n),
\mylabel{eq:etakdef}
\end{eqnarray}
which are linear combinations of the {\it $k$-th rank moments}
\begin{equation}
N_k(n) := \sum_{m = - \infty}^{\infty} m^k N(m,n)
\mylabel{eq:Nkdef}
\end{equation}
considered by
Atkin and the second author \cite{AG}.  Using the rank symmetry $N(-m,n) =
N(m,n)$, Andrews showed that $\eta_{2k+1}(n)=0$, and thus we need
only consider even rank moments. For these we define the rank
generating function
\begin{eqnarray*}
R_{k+1}(q):=
\sum_{n=0}^{\infty} \eta_{2k}(n)\,q^n.
\end{eqnarray*}

The function $R_2(q)$ was studied in detail by the first  author in
\cite{B1}.  One of the key results relates $R_2(q)$ to a certain weak
Maass form (see Section \ref{MaassSection}), but the connection is
more complicated than in the case of usual ranks due to double poles
in the generating function.  This leads to expressions involving {\it
quasimodular forms}, which are meromorphic functions $f:\H \to \C$
that can written as a linear combination of derivatives of modular
forms.
Furthermore, asymptotics and congruences for $\eta_2(n)$ are obtained
as applications of the modularity of the generating function.

In the present work we consider the case of general $k$. The
functions that arise in this setting require yet a more general
analytic definition; we say that $f:\H \to \C$ is a {\it quasimock
theta function}  if there exists a quasimodular form $h(q)$ such that
$f(q)+ h(q)$ is a linear combination of derivatives of the
holomorphic parts of weak  Maass forms.  Moreover we call linear
combinations of derivatives of weak Maass forms \textit{quasiweak
Maass form}.
\begin{theorem} \label{quasitheorem}
The function $q^{-1}R_{k+1}(q^{24})$ is  a quasimock theta function.
\end{theorem}
\begin{remark}
The highest weight component includes a weak Maass form of weight $2k - 1/2.$
\end{remark}

The idea of the proof of Theorem \ref{quasitheorem} is to relate
certain rank and crank moments via a differential equation (see
Section  \ref{RCSection}), and then argue inductively, using the fact
that the crank moment generating functions are quasimodular forms.
The base case $k=1$ is considered in \cite{B1}, although our
induction step actually requires a new ``twisted" version of those
results.

Theorem \ref{quasitheorem} has many applications. We only address
some of these here.
We first consider ccongruences for partition statistics. For this we
let $NF_k(r,t;n)$ be the
number of $k$-marked Durfee symbols of size $n$ with full rank
congruent to $r$ modulo $t$ (see Section \ref{CombResults}).
We show that the full rank satisfies infinitely many congruences,
just as the first author and Ono proved for Dyson's original rank
\cite{BO2}.
The case $k=2$ of the following theorem was proven in \cite{B1}.
\begin{theorem} \label{maintheorem}
Let $t$ be a positive odd integer, suppose that $j \in \N$, $k \geq
3,$ and let
$\mathcal{Q} \nmid 6t$  a prime. Then there exist infinitely many
arithmetic progressions $An+B$, such that for every $0 \leq r<t$, we
have
\begin{displaymath}
NF_k(r,t;An+B) \equiv 0 \pmod{\mathcal{Q}^j}.
\end{displaymath}
\end{theorem}

For the proof of Theorem \ref{maintheorem}, we  employ the fact that
the rank generating functions are holomorphic parts of weak Maass
forms, and also the conclusion from Theorem \ref{quasitheorem} for $R_k$.
Additional complications arise in our proof if $t$ has a prime
divisor $p_t$ that is small relative to $k$ (specifically, $p_t \leq
2k$).  To resolve this case, we extend a result from \cite{B1} and
prove the modularity of a certain ``twist" of the second moment
function
(see Sections \ref{RCSection} and  \ref{TwistSection}).

One nice consequence of Theorem \ref{maintheorem} is a combinatorial
decomposition of   congruences for $\eta_{2k}(n)$.
\begin{corollary}
Let $j \in \N$ and $Q>3$ a prime.
Then there exist infinitely many arithmetic progressions $An+B$ such
that
\begin{eqnarray*}
\eta_{2k}(An+B) \equiv 0 \pmod{  \mathcal{Q}^j }.
\end{eqnarray*}
\end{corollary}

\noindent To illustrate the nature of these arithmetic progressions, we give
some of the simpler examples
\begin{align*}
\eta_2(11^3 n + 479) &\equiv 0 \pmod{11},\\
\eta_4(11 n) &\equiv 0 \pmod{11},\\
\eta_6(49 n + 19) &\equiv 0 \pmod{7},\\
\eta_8(13^2 n + 162) &\equiv 0 \pmod{13}.
\end{align*}

See \cite{Gar} for a detailed explanation of the method used to find
these explicit congruences.  A significant part of the technique is
in showing that the rank moment generating functions are congruent
modulo $\mathcal{Q}$ to modular forms over a restricted set of
coefficients.  This is of independent interest, and is described
precisely in the following theorem.

\begin{theorem}\label{Congruencetheorem}
Suppose $\ell>3$ is prime.
Define $1\le \beta_{\ell}\le \ell-1$
such that $24\beta_{\ell}\equiv1\pmod{\ell}$ and let $r_{\ell} := \frac{24\beta_{\ell}-1}{\ell}.$
\begin{enumerate}
\item The generating function for the second rank moment satisfies
\beq
\sum_{n=0}^\infty N_2(\ell n + \beta_\ell) q^{24n+r_\ell}
\equiv \eta^{r_\ell}(24z) G_{\ell,2}(24z)
\pmod{\ell},
\notag
\eeq
where $G_{\ell,2}(z)$ is a sum of level $1$ modular forms with $\ell$-integral
coefficients, each of
weight at most $\frac{\ell(\ell+3)-r_\ell -1}{2}$.

\item
For $2\le k \le \tfrac{\ell-3}{2}$,
\beq
\sum_{n=0}^\infty N_{2k}(\ell n + \beta_\ell) q^{24n+r_\ell}
\equiv
c_k \sum_{n=0}^\infty N_2(\ell n + \beta_\ell) q^{24n+r_\ell}
+ \eta^{r_\ell}(24z) G_{\ell,2k}(24z)
\pmod{\ell},
\notag
\eeq
where $c_k$ is an integer, and $G_{\ell,2k}(z)$ is a sum of level $1$ modular forms with $\ell$-integral coefficients
and weight at most $k(\ell + 1) - 1 +\tfrac{1}{2}(\ell-r_\ell)$.

\item Finally,
\begin{align}
\sum_{n=0}^\infty N_{\ell-1}(\ell n + \beta_\ell) q^{24n+r_\ell}
&\equiv
c_{\ell-1} \sum_{n=0}^\infty N_2(\ell n + \beta_\ell) q^{24n+r_\ell}
+ \eta^{r_\ell}(24z) G_{\ell,\ell-1}(24z)
\notag \\
& \quad + \frac{1}{\ell} \eta^{r_\ell}(24z)
\left(
H_{1,\ell}(24z) - H_{2,\ell}(24z)
 \right)
\pmod{\ell},
\nonumber
\end{align}
where $G_{\ell,\ell-1}(z)$ is a sum of level $1$ integral modular forms
of weight
at most $ \frac{\ell(\ell+1)-r_\ell -3}{2}$;
$c_{\ell-1}$ is some integer; and
$H_{1, \ell}(z)$, $H_{2,\ell}(z)$ are integral modular forms of
weight $\tfrac{\ell(\ell-1)-r_\ell-1}{2}$ and
$\tfrac{\ell(\ell+1)-r_\ell-3}{2}$
respectively such that
\beq
H_{1,\ell}(z) \equiv H_{2,\ell}(z) \pmod{\ell}.
\notag
\mylabel{eq:Hcong}
\eeq

\end{enumerate}
\end{theorem}
We illustrate parts (1) and (2) of this theorem for the case $\ell=11$
(see \cite{Gar})
\begin{align}
\sum_{n=0}^\infty N_2(11n+6)q^{24n+13}
&\equiv 3 \eta^{13}(24z) \pmod{11}, \notag\\
\sum_{n=0}^\infty N_4(11n+6)q^{24n+13}
&\equiv 7 \eta^{13}(24z) \pmod{11}, \notag\\
\sum_{n=0}^\infty N_6(11n+6)q^{24n+13}
&\equiv \eta^{13}(24z)(4 + E_4(24z)) \pmod{11}, \notag\\
\sum_{n=0}^\infty N_8(11n+6)q^{24n+13}
&\equiv \eta^{13}(24z)(5 + 6 E_4(24z) + 6 E_6(24z)) \pmod{11}.
\notag
\end{align}

The theory of weak Maass forms can also be employed to show
identities for differences of rank moments.
Relations between non-holomorphic parts are responsible for the
existence of such identities.
To state our results, let
\begin{eqnarray*}
R_{r,s,t,d}^{(k)}(q) :=
\sum_{n=1}^{\infty}
\left(
NF_k(r,t;tn+d)
-NF_k(s,t;tn+d)
\right) q^{ 24(tn+d  )-1   }.
\end{eqnarray*}
We show that in certain cases this  function is  a weakly holomorphic
modular form whose poles (if there are any) are supported on the
cusps.
Similar results for Dyson's rank were shown in \cite{BOR}.
\begin{theorem} \label{IdentityTheorem}
Assume that $t \geq 5$ is a prime, $0 \leq r,s <t$, and $0 \leq
d<t$.
Then the following are true.
\begin{enumerate}
\item 
\label{Thm:R-1}
If $\leg{1-24d}{t}=-1$, then  $R_{r,s,t,d}^{(k)}(q)$ is a
quasimodular form  on
$\Gamma_1(576t^6)$.
If $k \leq \frac{p_t}{2}$, then it is  weakly holomorphic.
\item 
\label{Thm:R1}
If $\leg{1-24d}{t}=1$, and $2r,2s \not \equiv 3(1-d+2u)
\pmod{2t}$, $r,s \not \equiv 2+3u, 1+3u, 2-3d+3u,1-3d+3u \pmod t$ for
all $0 \leq u \leq d-1$, then  $R_{r,s,t,3d}^{(2)}(q)$ is a weakly
holomorphic modular form on
$\Gamma_1(576t^6)$.
\end{enumerate}
\end{theorem}
\noindent One can use Theorem \ref{IdentityTheorem} along with the valence
formula to prove concrete identities.
\begin{remark}
For each $t$ there exists at least one $0 \leq d<t$ such that the
statement of Theorem \ref{IdentityTheorem}  is nontrivial.
\end{remark}

The paper is organized as follows.
In Section \ref{CombResults} we recall facts about marked Durfee
symbols.
In Section \ref{MaassSection}, we give the connection between rank
generating functions and weak Maass forms. We then relate rank and
crank moments via a differential equation in Section
\ref{RCSection}.  Next, in Section \ref{TwistSection} we introduce a
certain twisted moment function that shows up in the case $k \geq
\frac{p_t}{2}$. Section \ref{MainSection}  is devoted the proofs of
Theorems \ref{quasitheorem} and  \ref{maintheorem}.
In Section \ref{CongruenceSection} we prove Theorem  \ref{Congruencetheorem}
by first making a detailed $\ell$-adic analysis of  the
rank-crank moment relation.
Section  \ref{IdentSection} provides identities for rank
differences.
\section{Combinatorial results on marked Durfee symbols}
\label{CombResults}
Here we give some of the facts about marked Durfee symbols shown in
\cite{An2}.
Recall that the largest square of nodes in the Ferrers graph of a
partition is called the \textit{Durfee square}.
The \textit{Durfee symbol} consists of 2 rows and a subscript, where
the top row consists of the columns to the right of the Durfee
square, the bottom row consists of the rows  below the Durfee square
and the subscript denotes the side length of the Durfee square.
The number being partitioned is equal  to the sum of the rows of the
symbol plus the number of nodes in the Durfee square.  Note that the
parts in both rows must be non-increasing.  As an example, the Durfee
symbol
\begin{eqnarray*}
\left(
\begin{matrix}
2& \\
3&3&1
\end{matrix}
\right)_4
\end{eqnarray*}
represents a partition of $2+3+3+1+4^2=25$.

Andrews defined \textit{k-marked Durfee symbols} by using $k$
distinct copies (or colors) of the integers designated by $\{
1_1,2_1,\cdots \}$, $\{ 1_2,2_2,\cdots \},\cdots$  $,\{
1_k,2_k,\cdots \}$.  We form Durfee symbols as before and use the $k$
copies of integers for parts in both rows.
We additionally demand that:
\begin{enumerate}
\item
The sequence of subscripts in each row are non-increasing.
\item Each of the subscript $1, \cdots, k-1$ occurs at least once in
the top row.
\item If $M_1, \cdots, M_{k-1}$ are the largest parts with  their
respective subscripts in the top row,  then all parts in the bottom
row with subscript 1 lie in $[1,M_1]$,   with subscript $2$ lie in
$[M_1,M_2], \cdots$,
and with subscript $k$ lie in $[M_{k-1},S]$,
where $S$ is the side of the Durfee square.
\end{enumerate}
The {\it size} of a $k$-marked Durfee symbol is simply the size of
the partition that is obtained by ignoring the colors; we let
$\mathcal{D}_k(n)$ denote the number of $k$-marked Durfee symbols of
size $n$.

In \cite{An2} Andrews showed that $k$-marked Durfee symbols arise
naturally in the combinatorial study of the rank moment functions; in
particular, for $k \geq 1$,
\begin{eqnarray*}
\mathcal{D}_{k+1} (n) = \eta_{2k}(n).
\end{eqnarray*}
He also proved some striking congruences for Durfee symbols,
including
\begin{align*}
\calD_2(5n + a) & \equiv 0 \pmod{5} \qquad a \in \{1, 4\}, \\
\calD_2(7n + a) & \equiv 0 \pmod{5} \qquad a \in \{1, 5\}, \\
\calD_3(7n+a) & \equiv 0 \pmod{7} \qquad a \in \{1, 5\}.
\end{align*}

Following Dyson's lead, Andrews next associated a collection of
``ranks" to $k$-marked Durfee symbols. For such a Durfee symbol
$\delta$, we define the \textit{full rank} $FR(\delta)$ by
\begin{eqnarray*}
FR(\delta) := \rho_1(\delta) + \cdots + k \rho_k(\delta),
\end{eqnarray*}
where the $i$-th rank $\rho_i(\delta)$ is given by
\begin{eqnarray*}
\rho_i(\delta) := \left\{
\begin{array}{ll}
\tau_i(\delta) - \beta_i(\delta)-1 & \text{ for }1 \leq i<k,\\
\tau_i(\delta)- \beta_i(\delta) &\text{ for } i=k.
\end{array}
\right.
\end{eqnarray*}
Here $\tau_i(\delta)$ (resp. $\beta_i(\delta)$) denotes the number of
entries in the top (resp. bottom) row of $\delta$ with subscript
$i$.
We let $NF_k(m,n)$ denote the number of k-marked Durfee symbols of
size $n$
with full rank $m$, and $NF_k(r,t;n)$ denote the number of $k$-marked
Durfee symbols of size $n$ with full rank congruent to $r$  modulo
$t$.
Finally, define the generating function
\begin{eqnarray*}
R_k(w;q):=
\sum_{n=1}^{\infty} \sum_{m \in \Z} NF_k(m,n)\, w^m\, q^n.
\end{eqnarray*}
In particular
\begin{eqnarray*}
R_k(1;q)=R_k(q).
\end{eqnarray*}

\section{Ranks and weak Maass forms} \label{MaassSection}

Here we recall results of \cite{BO2} and \cite{B1}.
Let us first give  the
definition of a weak Maass form.
If $k\in \frac{1}{2}\Z\setminus
\Z$,
$z=x+iy$ with $x, y\in \R$, then the weight $k$ hyperbolic
Laplacian is given by
\begin{equation}\mylabel{eq:laplacian}
\Delta_k := -y^2\left( \frac{\partial^2}{\partial x^2} +
\frac{\partial^2}{\partial y^2}\right) + iky\left(
\frac{\partial}{\partial x}+i \frac{\partial}{\partial y}\right).
\end{equation}
If $v$ is odd, then define $\epsilon_v$ by
\begin{equation}
\epsilon_v:=\begin{cases} 1 \ \ \ \ &{\text {\rm if}}\ v\equiv
1\pmod 4,\\
i \ \ \ \ &{\text {\rm if}}\ v\equiv 3\pmod 4. \end{cases}
\end{equation}
Moreover, we let $\chi$ be a Dirichlet character.

A {\it (harmonic) weak Maass form of weight $k$ with Nebentypus
$\chi$ on a subgroup
$\Gamma \subset \Gamma_0(4)$} is any smooth function $g:\H\to \C$
satisfying the following:
\begin{enumerate}
\item For all $A= \left(\begin{smallmatrix}a&b\\c&d
\end{smallmatrix} \right)\in \Gamma$ and all $z\in \H$, we
have
\begin{displaymath}
g(Az)= \leg{c}{d}^{2k}\epsilon_d^{-2k} \chi(d)\,(cz+d)^{k}\ g(z).
\end{displaymath}
\item We  have that $\Delta_k g=0$.
\item The function $g(z)$ has
at most linear exponential growth at all the cusps of $\Gamma$.
\end{enumerate}
Now let $0 <a<c$ and define
\begin{eqnarray*}
D\left( \frac{a}{c};z\right) &:=& - S \left( \frac{a}{c};z\right)
+q^{-\frac{\ell_c}{24} }\, R \left(\zeta_c^a;q^{ \ell_c} \right),\\
S\left( \frac{a}{c};z\right) &:=& - \frac{i \sin \left( \frac{\pi
a}{c}\right) \ell_c^{\frac{1}{2} }}{\sqrt{3}}\,
\int_{-\bar z}^{i \infty}
\frac{\Theta\left( \frac{a}{c};\ell_c \tau\right)}{ \sqrt{ i (\tau+z)
} } \, d \tau,
\end{eqnarray*}
where  $\ell_c:= \lcm (2c^2,24)$, $\zeta_c:= e^{\frac{2 \pi i}{c} }$,
and  $\Theta\left(\frac{a}{c};\tau \right)$ is a certain weight
$\frac{3}{2}$ cuspidal theta function (for the exact definition  see
\cite{BO2}).

\begin{theorem} \label{MaassformTheorem}
If $0 <a<c$, then $D \left(\frac{a}{c};z \right)$  is a weak Maass
form of weight $\frac12$ on
$\Gamma_c:=  \left<\left(\begin{smallmatrix} 1 &1\\0 &1
\end{smallmatrix}\right),  \left( \begin{smallmatrix} 1&0\\\ell_c^2&1
\end{smallmatrix} \right) \right> $.
If $c$ is odd, then it is on $\Gamma_1\left(6f_c^2l_c \right)$, where
$f_c:= \frac{2c}{\gcd(c,6)}$.
Its non-holomorphic part has the expansion
\begin{eqnarray*}
- \frac{2}{\sqrt{\pi}} \sin \left(\frac{\pi a}{c} \right)
\sum_{m \pmod{f_c }} (-1)^m \sin \left( \frac{\pi a(6m+1)}{c}\right)
\sum_{ n \equiv 6m+1 \pmod{6f_c}} \Gamma \left(\frac12;\frac{\ell_c
n^2y}{6} \right)\, q^{ -\frac{\ell_cn^2}{24}},
\end{eqnarray*}
where
$$
\Gamma (\alpha;x):= \int_{x}^{\infty}e^{-t}\, t^{a-1}\, dt.
$$
\end{theorem}
We next turn to $R_2(q)$.
Define
\begin{eqnarray*}
\mathcal{R}(z) := R_2 (24 z) e^{ - 2\pi i z},
\end{eqnarray*}
\begin{eqnarray}\mylabel{eq:Mordelint}
\mathcal{N}(z):=
\frac{i}{4\sqrt{2}\pi}
\int_{- \bar z}^{i \infty}
\frac{\eta(24\tau )}{(- i (\tau+z))^{\frac32}} \, d\tau,
\end{eqnarray}
and
\begin{eqnarray}
\mathcal{M} (z) :=
\mathcal{R}(z)-
\mathcal{N}(z)  -
\frac{1}{24\eta(24 z)}
+
\frac{E_2(24 z)}{8\eta(24 z)},
\end{eqnarray}
where  as usual
$$E_2(z):=
1- 24 \sum_{n=1}^{\infty} \sigma_1(n)\,q^n
$$
with $\sigma_1(n) := \sum_{d|n}d$.
This function is a quasimodular form.

The following result is shown in \cite{B1}.
\begin{theorem} \label{maintheorem2}
The function $\mathcal{M}(z)$ is a harmonic weak Maass form of weight
$\frac{3}{2}$
on  $\Gamma_0(576)$ with Nebentypus character
$\chi_{12}:=\leg{12}{\cdot}$.
Its non-holomorphic part has the expansion
$$
\mathcal{N}(z) = \frac{1}{4 \sqrt{\pi}} \sum_{k \in \Z} (-1)^k (6k+1)
\Gamma \left(-\frac12;4 \pi (6k+1)^2y \right)\, q^{ -(6k+1)^2}.
$$
\end{theorem}

\section{Relation between rank  and crank moments} \label{RCSection}
In this section we recall certain relations between
rank  and crank moments and consider twisted generalisations. For
details we refer the reader to \cite{AG}.
The $j$-th rank moment $N_{j}$ is defined in \eqn{Nkdef}.
Note that the symmetrized moment $\eta_{2k}(n)$ can easily be written
as a linear combination of $N_{2j}(n)$ with $j \leq k$ (again,
$N_{2j+1}(n) = 0$ due to symmetry).
Define the generating function
\begin{eqnarray*}
\mathcal{R}_j(q) &:=& \sum_{n \geq 1} N_j(n) q^n.
\end{eqnarray*}
We will see that these functions are related to certain
quasiweak Maass forms.

We next consider crank moments.
 Recall that   the \textit{crank} of a partition is defined to be the
 largest part if the partition contains no ones, and is otherwise the
 difference between the number of parts larger than the number of
 ones and the number of ones.
 For $n >1$, we denote by $M(m,n)$ the number of partitions of $n$
 with crank equal to $m$, and define the boundary values by $M(0,1)
 :=-1$, $M(-1,1)
 :=M(1,1):=1$, with $M(m,1):=0$ otherwise.
 The generating function for the crank is then
 \begin{eqnarray*}
 C(w;q):=
 \sum_{n \geq 0} \sum_{m \in \Z} M(m,n)\,w^m\, q^n
 = \prod_{n=1}^{\infty}
 \frac{\left( 1-q^n\right)}{\left( 1-wq^n\right)\left(1-w^{-1}q^n
 \right)} .
 \end{eqnarray*}
 This function is essentially a modular form when $w$ is the root of
 unity $\zeta^a_c$.  The numerator is $\eta(z) q^{-1/24}$, and the
 denominator is an algebraic integer times a weight zero {\it Siegel
 function} of level $2c^2$ (see \cite{KL}); this implies that $q^{-1}
 C \left( \zeta_c^a;q^{24}\right)$ is a weight $\frac12$ weakly
 holomorphic modular form on $\Gamma_1(2 \cdot \lcm(c^2,288))$.

Analogous to the development above, define the  \textit{$j$-th crank
moment}  as
\begin{eqnarray*}
M_j(n) := \sum_{k \in \Z} k^j M(k,n),
\end{eqnarray*}
which again satisfies $M_{2j+1}(n) =0$.  Denote the crank moment
generating function by
\begin{eqnarray*}
C_j(q) &:=& \sum_{n \geq 1} M_j(n) q^n.
\end{eqnarray*}

Now define the differential operators
\begin{eqnarray*}
\delta_{q} = \delta_{z} &:=& q \: \frac{d}{dq} = \frac{1}{2 \pi i }
\frac{d}{d z},\\
\delta_w&:=& w \frac{d}{dw}.
\end{eqnarray*}
In \cite{AG}, Atkin and the second author derived a recurrence relation for the
functions $C_a$:
\begin{eqnarray}\mylabel{eq:relation}
C_a(z) = 2 \sum_{ j=1}^{\frac{a}{2}-1 }
\left(
\begin{matrix}  a-1\\2j-1 \end{matrix}
\right)
\Phi_{2j-1 } (z) C_{a-2j}(z)
+ 2 \Phi_{a-1}(z) P(z).
\end{eqnarray}
Here
\begin{eqnarray*}
\Phi_{2j-1} (z):=
\sum_{n=1}^{\infty} \sigma_{2j-1}(n)q^n,
\end{eqnarray*}
where
$\sigma_j(n):=\sum_{d|n} d^j$.
These functions are simply a rescaling of the classical weight $j$
Eisenstein series minus their constant terms, since
\begin{eqnarray*}
E_j(z):=
1- \frac{2j}{B_j} \Phi_{j-1}(z),
\end{eqnarray*}
where $B_j$ is the $j$th Bernoulli number.
For even $j>2$, the function $E_j(z)$ is a  modular form of level
$1$, whereas $E_2(z)$ is a quasimodular form.
Thus we conclude inductively from (\ref{eq:relation}) that
$q^{-1}C_a(q^{24})$ is a quasimodular form.

Atkin and the second author also proved a differential equation for the crank and rank generating
functions, called the ``rank-crank PDE":
\begin{equation}
\label{eq:PDE}
w(q;q)_\infty^2 C(w;q)^3=
\left(3(1-w)^2 \delta_q + \frac{1}{2}(1-w)^2 \delta_w^2 - \frac{1}{2}(w^2-1)\delta_w + w \right)R(w;q).
\end{equation}
For
our current purposes, we are most interested in an identity that they derived
by repeatedly applying $\delta_w$ to (\ref{eq:PDE}) and setting $w=1$, namely that for $a \geq 2,$
\begin{multline}\mylabel{eq:rankcrank}
\sum_{i=0}^{a/2-1}
\left(
\begin{matrix}
a\\2i
\end{matrix}
\right)
\sum_{ \substack{\alpha+ \beta + \gamma = a-2i   \\ \alpha, \beta,
\gamma \geq 0 \text{ even} }}
\left(
\begin{matrix}
a-2i\\\alpha,\beta, \gamma
\end{matrix}
\right)
C_{\alpha}(z) C_{\beta} (z) C_{\gamma}(z) P^{-2}(z)  -3(2^{a-1 } -1  ) C_{2}(z)  \\
= \frac12(a-1)(a-2) \mathcal{R}_{a}  (z)
+ 6 \sum_{i=1}^{a/2-1 }
\left( \begin{matrix} a\\ 2i   \end{matrix} \right)
( 2^{2i-1 }-1) \delta_q\left( \mathcal{R}_{a-2i}  (z)\right)\\
+ \sum_{i=1}^{ a/2 -1} \left[
\left(  \begin{matrix} a\\ 2i+2  \end{matrix} \right)
(2^{2i+1 } -1)
-2^{ 2i}
\left(  \begin{matrix} a\\ 2i+1  \end{matrix} \right)
+  \left(  \begin{matrix} a\\ 2i \end{matrix} \right)
\right] \mathcal{R}_{a-2i}(z).
\end{multline}
Therefore modularity properties of $\mathcal{R}_a$ can be inductively
concluded from modularity properties of $C_a$ and $\mathcal{R}_2$.
In particular, the functions $\mathcal{R}_{a}$ are related to
the quasiweak Maass forms that we desribe in more detail later.

We next turn to twisted rank and crank moments.
We  have
\begin{eqnarray*}
\delta_w \left(C(w;q) \right)= L(w;q)\, C(w;q),
\end{eqnarray*}\
where
\begin{eqnarray*}
L(w;q):=
\sum_{n,m \geq 1}
\left(
w^m q^{nm}
-w^{-m} q^{nm}
\right).
\end{eqnarray*}
Thus
\begin{eqnarray}\mylabel{eq:deltaL}
\left[ \delta^j_w L(w;q) \right]_{w=\zeta}
= \sum_{n,m \geq 1}
\left(
m^j \zeta^m q^{nm}
- (-m)^j \zeta^{-m}q^{nm}
\right).
\end{eqnarray}

Now recall the theory of Eisenstein series on congruence subgroups
(see section III.3 in \cite{Kob}).  The $(0,a)$ Eisenstein series of
weight $j+1$ and level $c$ is given by
\begin{eqnarray*}
G_{j+1}^{(0,a)}(z) := b_{j+1}^{(0,a)} + c_{j+1}\sum_{n \geq 1}
\left( \sum_{d \mid n} d^j\left(\zeta_c^{ad} - (-1)^j
\zeta_c^{-ad}\right)\right)q^n,
\end{eqnarray*}
where
\begin{eqnarray*}
b_{j+1}^{(0,a)} = \sum_{\begin{subarray}{c} n \geq 1 \\ n \equiv a
\pmod{c} \end{subarray}} n^{-j-1} +
\sum_{\begin{subarray}{c} n \geq 1 \\ n \equiv -a
\pmod{c}\end{subarray}} (-n)^{-j-1}\cdot
\end{eqnarray*}
and
\begin{eqnarray*}
c_{j+1} = \frac{2(j+1)(-1)^{j}\zeta(j+1)}{c^{j+1}B_{j+1}}.
\end{eqnarray*}
These Eisenstein series are in $M_{j+1}(\Gamma_1(c))$ for $j \geq 2$ (as before, they may be quasimodular
at weight $2$), and thus the
series from (\ref{eq:deltaL}) is again a rescaled quasimodular form minus its
constant coefficient.  Note also that the constant
$b_{j+1}^{(0,a)}/c_{j+1}$ must be an algebraic integer since all of
the other terms in the rescaled series are.

Using (\ref{eq:rankcrank}) and results from Section
\ref{TwistSection}, one can prove
modularity properties of the {\it twisted moment functions}
\begin{eqnarray*}
\mathcal{R}_{j,a,c}(q)&:=& \sum_{n \in \N} \mathcal{R}_j \left(
\frac{a}{c};n\right)\, q^n,\\
C_{j,a,c}(q)&:=& \sum_{n \in \N} C_j \left( \frac{a}{c};n\right)\, q^n,
\end{eqnarray*}
where
\begin{eqnarray*}
\mathcal{R}_j \left( \frac{a}{c} ; n\right)&:= &\sum_{k \in \Z}
k^j\zeta_c^{ak} N(k,n),\\
C_j \left( \frac{a}{c} ; n\right)&:= &\sum_{k \in \Z} k^j\zeta_c^{ak}
M(k,n).
\end{eqnarray*}

\section{A twisted moment function} \label{TwistSection}
Define for coprime integers $0<a<c$ the twisted second  moment rank
generating function
\begin{eqnarray}\mylabel{eq:TwistGen}
R_2 \left(\frac{a}{c};q \right)
:= \frac{\zeta_{2c}^a}{2(q;q)_{\infty}}
\sum_{ n \in \Z}
\frac{(-1)^{n+1} q^{ \frac{n}{2}(3n+1)   }}{\left( 1-\zeta_c^aq^n
\right)}
+
\frac{\zeta_{2c}^{3a}}{(q;q)_{\infty}}
\sum_{ n \in \Z}
\frac{(-1)^{n+1} q^{ \frac{3n}{2}(n+1)   }}{\left( 1-\zeta_c^aq^n
\right)^2}.
\end{eqnarray}
We relate this function to a weak Maass form.
For this, let
\begin{eqnarray*}
\mathcal{M}_{\frac{a}{c}} (z)
= \mathcal{R}\left(  \frac{a}{c};q  \right) +
\mathcal{N}_{\frac{a}{c}} (z).
\end{eqnarray*}
Here
\begin{eqnarray*}
\mathcal{R} \left(\frac{a}{c};q \right)&:=& q^{-1} R_2
\left(\frac{a}{c};q^{24} \right),\\
\mathcal{N}_{\frac{a}{c}}(z)
&:=& -\frac{i}{64\sqrt{3\pi} }
\int_{-\bar z}^{i \infty}
\frac{ \Theta_{ a,c } \left( -\frac{1}{16d_c^2 \tau }\right)(-i\tau
)^{ -\frac12}}{(  -i (   \tau + z  )   )^{ \frac32}} \, d \tau,
\end{eqnarray*}
where  $d_c:= \lcm (6,c)$, and where
\begin{eqnarray*}
\Theta_ {a,c}(\tau):= \sum_{m \equiv \frac{d_c}{6}  \pm \frac{d_c
a}{c} \pmod{d_c }}
(-1)^m \, e^{ 2 \pi i m^2 \tau}.
\end{eqnarray*}
\begin{theorem} \label{twistTheorem}
The function $\mathcal{M}_{\frac{a}{c}}(z) $ is a weight $\frac32$
weak Maass form on $\Gamma_1\left(96d_c^2\right)$.
\end{theorem}
The first step in proving Theorem \ref{twistTheorem} is to show a
transformation law for \
$R_2 \left(\frac{a}{c};q \right)$. Due to double poles, we cannot
work directly with this function, but use a function of an additional
parameter $w$ that is related and only has single poles.
Define
\begin{eqnarray*}
R_2 \left(\frac{a}{c},q;w  \right)
:= \zeta_{2c}^a \frac{e^{\pi i w }}{(q;q)_{\infty}}
\sum_{n  \in \Z}
\frac{(-1)^{ n+1 } q^{ \frac{n}{2}(3n+1) }}{\left(  1-\zeta_c^ae^{2
\pi iw}  q^n  \right)}.
\end{eqnarray*}
This function is connected to $R_2 \left(\frac{a}{c};q \right)$ by
\begin{eqnarray*}
L \left( R_2 \left( \frac{a}{c},q;w   \right) \right)
= R_2 \left(\frac{a}{c};q\right),
\end{eqnarray*}
where for  a function $g$ that is differentiable  in some
neighborhood of $0$,  we define
\begin{eqnarray*}
L(g):=\left[  \frac{1}{2\pi i} \frac{\partial}{\partial w}  g(w)
\right]_{w=0}.
\end{eqnarray*}
We prove a transformation law for  $R_2 \left(\frac{a}{c},q;w\right)$
and then apply $L$.
Some of the calculations are  similar to those in \cite{B1},
therefore we skip some of the details here and refer the reader to
that paper.

We begin with some useful notation.
For $\frac{a}{c} \not \in \left\{ 0,\frac16,\frac12,\frac56
\right\} $, let
\begin{eqnarray*}
s=s(a,c):=
\left\{
\begin{array}{ll}
0&\text{if } 0 < \frac{a}{c} <\frac16,\\[0.3ex]
1&\text{if } \frac16< \frac{a}{c} <\frac12,\\[0.3ex]
2&\text{if }  \frac12< \frac{a}{c} <\frac56,\\[0.3ex]
3&\text{if } \frac56<\frac{a}{c}<1.
\end{array}
\right.
\end{eqnarray*}
Moreover define
\begin{eqnarray*}
\omega_{h,k} :=
\exp\left(\pi i t(h,k) \right),
\end{eqnarray*}
where we have denoted the standard Dedekind sum by
\begin{eqnarray*}
t(h,k):= \sum_{\mu \pmod k}  \left( \left( \frac{\mu}{k}\right)
\right)   \left( \left( \frac{h \mu}{k}\right) \right),
\end{eqnarray*}
with
\begin{eqnarray*}
((x)):= \left \{
\begin{array}{ll}
x- \lfloor x \rfloor - \frac{1}{2} &\text{if } x \in \R \setminus \Z
,\\
0&\text{if } x \in \Z.
\end{array}
\right.
\end{eqnarray*}
Let
\begin{eqnarray*}
T_2 \left(\frac{a}{c},q;w \right)
&:=& \frac{1}{(q;q)_{\infty}}
\sum_{\pm} \pm
e^{ \pm \pi i w}\, q^{\pm \frac{a}{2c}} \sum_{ m \geq 0} (-1)^m
\frac{ q^{ \frac{m}{2} (3m+1) \pm ms} }{1-e^{\pm 2 \pi i w } q^{\pm
\frac{a}{c}+m}} ,\\
I_{k, \nu,a,c}^{\pm}  (z;w)&:= &\pm \int_{\R}
\frac{e^{ -\frac{3 \pi zx^2}{k}}}{\sinh \left( \frac{\pi zx}{k}
+\frac{\pi i}{6k} -  \frac{\pi i\nu}{k} \mp \pi i \left(  w+
\frac{a}{c}\right)\right)}\, dx.
\end{eqnarray*}
\begin{theorem}  \label{twistTransTheorem}
Assume the notation above.
Moreover for coprime integers $h$ and $k$, with $k>0$ and  either
$k=1$ or $2c^2|k$, let
$q:=e^{\frac{2 \pi i}{k}(h+iz)}$ and $q_1:=e^{\frac{2 \pi
i}{k}\left(h' + \frac{i}{z} \right)}$, with $z \in \C$,
$\text{Re}(z)>0$, where $h'=0$ for $k=1$, and $hh' \equiv -1
\pmod{2k}$ and  $h' \equiv 1 \pmod{2c^2}$ for
$2c^2|k$.
Then
\begin{eqnarray*}
R_2 \left( \frac{a}{c},q;w \right)
= S_1  -\frac{  z^{    \frac12}}{2k}\omega_{h,k}\,e^{-\frac{\pi
z}{12k} }
\sum_{  \substack{  \nu \pmod k \\   \pm}  } (-1)^{\nu}
e^{ \frac{\pi i h'}{k}( -3 \nu^2 +\nu)}  I_{k,\nu,a,c}^{\pm}(z;w),
\end{eqnarray*}
where
\begin{displaymath}
S_1:=
\left\{
\begin{array}{ll}
- i z^{-\frac12}   e^{\frac{\pi}{12} \left(z^{-1}-z \right) - \frac{2
\pi s}{z}\left( w+\frac{a}{c} \right) +\frac{3 \pi}{z} \left(
w+\frac{a}{c}\right)^2 }
T_2 \left(\frac{a}{c},q_1;\frac{iw}{z} \right)
&\text{if } k=1,\\
-\frac{i}{z^{\frac{1}{2}}  } \, \omega_{h,k} e^{ \frac{3 \pi
kw^2}{z}  +\frac{\pi}{12k}\left( z^{-1} -z\right) }
R_2 \left(\frac{a}{c},q_1;\frac{w}{iz}  \right)&\text{if } 2c^2 |k.
\end{array}
\right.
\end{displaymath}
\end{theorem}
\begin{proof}
Let
\begin{eqnarray*}
\widetilde{R}_2
\left( \frac{a}{c},q; w\right)
:= \zeta_{2c}^a\, e^{ \pi iw }
\sum_{n  \in \Z}
\frac{(-1)^{ n+1 } q^{ \frac{n}{2}(3n+1) }}{\left(  1-\zeta_c^ae^{2
\pi iw}  q^n  \right)}.
\end{eqnarray*}
Poisson summation yields
\begin{equation}\mylabel{eq:poissonsum}
\widetilde{R}_2 \left( \frac{a}{c},q; w\right)
= -\frac{1}{2k}
\sum_{  \substack{  \nu \pmod{k}  \\  \pm  }  }
\pm (-1)^{\nu} e^{\frac{3 \pi i h\nu^2}{k} }
\sum_{  n \geq 0 }
\int_{ \R}
\frac{e^{   -\frac{3\pi zx^2}{k}   +\frac{\pi i}{k}(2n+1)(x - \nu)
}}{\sinh  \left(\frac{\pi z x}{k}  -\frac{\pi i h \nu}{k}  \mp \pi i
\left( w+\frac{a}{c}\right) \right)}\, dx.
\end{equation}
We shift the path of integration through
$\omega_n:=\frac{(2n+1)i}{6z}$. Using the residue theorem yields
\begin{eqnarray*}
\widetilde{\mathcal{R}}_2 \left(\frac{a}{c},q;w \right)
= \sum_1 + \sum_2.
\end{eqnarray*}
Here
\begin{eqnarray*}
\sum_1:= 2 \pi i \sum_{\text{residues}} ,
\end{eqnarray*}
where the sum runs over all residues of the integrand in
(\ref{eq:poissonsum}), and  where
\begin{eqnarray*}
\sum_2&:=&  -\frac{1}{2k}
\sum_{  \substack{  \nu \pmod{k}  \\  \pm  }  }
\pm (-1)^{\nu} e^{\frac{3 \pi i h\nu^2}{k} }
\sum_{  n \geq 0 }
\int_{\R  +\omega_n}
\frac{e^{   -\frac{3\pi zx^2}{k}   +\frac{\pi i}{k}(2n+1)(x - \nu)
}}{\sinh  \left(\frac{\pi z x}{k}  -\frac{\pi i h \nu}{k}  \mp \pi  i
\left( w+\frac{a}{c}\right) \right)}\, dx.
\end{eqnarray*}

To compute $\sum_2$, we observe that for  the computation of $\sum_2$
in \cite{B1} one does not need the fact that  $w$ is  small,
therefore we may change
$w \mapsto w +\frac{a}{c}$.
This yields
\begin{eqnarray*}
\sum_2 =
-\frac{  (q_1;q_1)_{ \infty}}{2k} \, e^{ -\frac{\pi}{12kz}  }
\sum_{\substack{  \nu \pmod k \\ \pm  } } (-1)^{ \nu}\, e^{ \frac{\pi
i h'}{k}( -3 \nu^2 +\nu)}
I_{ k,\nu,a,c}^{\pm}(z;w).
\end{eqnarray*}

We next turn to $\sum_1$.
First we consider the case  $k=1$.
In this case
poles of the integrand can only lie in points
$$
x_m^{\pm} :=\frac{i}{z} \left(m \pm \left(\frac{a}{c} +w \right)
\right).
$$
If we shift the path of integration through $\omega_n$, we have to
take those $x_m^{\pm}$ into account
for which $n \geq 3m \pm s$ and $m \geq \frac12(1 \mp1)$.
We denote the residues of each summand by $\lambda_{n,m}^{\pm}$.
Then one can easily see that
\begin{eqnarray*}
\lambda_{n,m}^{\pm}= \mp
\frac{e^{ - 3 \pi zx_m^{\pm 2} + \pi i (2n+1)x_m^{\pm}}}{2 \pi z\,
\cosh\left(\pi z x_m^{\pm} \mp  \pi i \left( w + \frac{a}{c}
\right)\right)}.
\end{eqnarray*}
Using that $\lambda_{n+1,m}^{\pm}=e^{ 2 \pi i x_m^{\pm}}
\lambda_{n,m}^{\pm}$,  one can compute that
\begin{multline*}
\sum_1 = \frac{1}{iz}
\, e^{\frac{3 \pi }{z} \left(w+\frac{a}{c} \right)^2 - \frac{2 \pi
s}{z} \left( w+ \frac{a}{c}\right)  }
\left(
\sum_{m \geq 0}  (-1)^m
\frac{e^{ - \frac{3 \pi m^2}{z} - \frac{2 \pi sm}{z} -\frac{\pi m}{z}
-\frac{\pi}{z}  \left( w + \frac{a}{c}\right)  }}{1-e^{ -\frac{2
\pi}{z}\left( w+\frac{a}{c}\right) - \frac{2 \pi m}{z} }} \right. \\
\left. -
\sum_{m > 0} (-1)^m
\frac{e^{ - \frac{3 \pi m^2}{z} +\frac{2 \pi sm}{z} -\frac{\pi m}{z}+
\frac{\pi}{z}  \left( w + \frac{a}{c}\right)  }}{1-e^{ \frac{2
\pi}{z}\left( w+\frac{a}{c}\right) - \frac{2 \pi m}{z} }}
\right).
\end{multline*}

We next consider the case $2c^2|k$.
Define the entire function
\begin{eqnarray*}
S_w^{\pm} (x):=
\frac{\sinh (x \pm \pi i kw) }{ \sinh \left(\frac{x}{k} \pm \pi i w
\right)  }.
\end{eqnarray*}
From this  one   can see that  poles of the integrand in
(\ref{eq:poissonsum}) only lie  in points
$$
x_m^{\pm}:=\frac{i}{z} (m \pm kw  ),
$$
and a non-trivial residue occurs for at most one $\nu$ modulo $k$,
which we may chose as
$$
\nu_m^{\pm} := -h' \left(m \mp \frac{ak}{c} \right).
$$
Shifting  the path of integration through $\omega_n$, we have  to
take those $m$ into account   for which $n \geq 3m \geq \frac{1}{2}
(1 \mp 1)$.
A lengthy calculation using the same methods as before gives
\begin{eqnarray*}
\sum_1:=
\frac{1}{iz} \, \zeta_{2c}^a\, e^{ \frac{3 \pi k w^2}{z}+\frac{\pi
w}{z}}
\sum_{ m \in \Z}
\frac{(-1)^{m+1}\, q_1^{\frac{m}{2} (3m+1) } }{1-\zeta_c^a\, e^{
\frac{2 \pi w}{z}}  q_1^m}.
\end{eqnarray*}
Now the theorem follows using
\begin{eqnarray*}
(q_1;q_1)_{\infty}
= \omega_{h,k} z^{\frac12} e^{ \frac{\pi}{12k} \left(z^{-1}-z
\right)}\, (q;q)_{\infty}.
\end{eqnarray*}
\end{proof}

Next we   realize the integrals occurring in Theorem
\ref{twistTransTheorem}  as theta integrals.
For this let
\begin{eqnarray*}
I_{ a,c}^{\pm} (w;z):=
\zeta_{2c}^a e^{ \pi i w}e^{ \mp\frac{\pi i}{6} } \int_{ \R}
e^{-\frac{3 \pi i x^2}{z} }
\frac{e^{ \mp \frac{\pi i x}{z}   }}{1-e^{  \mp\frac{\pi i}{3}  + 2
\pi i \left(  w +\frac{a}{c} \right) \mp \frac{2 \pi i x}{z} } }\,
dx.
\end{eqnarray*}
As in \cite{B1}, we show.
\begin{lemma} \label{thetaLemma}
We have
\begin{eqnarray*}
L \left( I_{a,c}^+(w;z) + I_{a,c}^-(w;z)\right)
= \frac{\sqrt{3}(-iz)^2}{4\pi}  \int_{0}^{\infty}
\frac{\Theta_{a,c}\left(\frac{3 i u}{2 d_c^2}\right)}{( -i (iu + z) )
^{\frac32}} \, du.
\end{eqnarray*}
\end{lemma}
To finish the proof of Theorem   \ref{twistTheorem},
we change in Theorem \ref{twistTransTheorem}  $z \mapsto \frac{i}{z}$
and   apply the operator $L$ on both sides.
Observe that
\begin{eqnarray*}
L \left(\frac{e^{ 2 \pi i s \left(w + \frac{a}{c} \right) z - 3 \pi i
z \left(w+\frac{a}{c} \right)^2 \pm \pi i z w  }}{1 - e^{\pm 2 \pi i
w z}q_1^{\pm \frac{a}{c} +m }}  \right)
= z \, e^{2 \pi i s \frac{a}{c} z- 3 \pi i \frac{a}{c} z }
\left(
\frac{\left( s- \frac{3a}{c} \pm \frac12 \right)}{1-q_1^{ m  \pm
\frac{a}{c}}}
\pm \frac{q_1^{m \pm \frac{a}{c} }}{\left(1-q_1^{m \pm \frac{a}{c} }
\right)^2}
\right).
\end{eqnarray*}
In the sum over the $-$ sign, we change $m \mapsto -m$, yielding
\begin{eqnarray*}
L \left(
e^{ 2 \pi i s \left(w + \frac{a}{c} \right)z - 3 \pi i z
\left(w+\frac{a}{c} \right)^2} \, T_2 \left(\frac{a}{c},q ; zw
\right) \right)
= z \, q^{ \frac{sa}{c} -\frac{3 a^2}{2c^2} } T_2\left(
\frac{a}{c};q\right),
\end{eqnarray*}
where
\begin{eqnarray*}
T_2 \left( \frac{a}{c};q\right)
:= \frac{1}{(q;q)_{\infty} }
\left( \left(s-\frac{3a}{c} \right)
q^{ \frac{a}{2c}}
\sum_{m \in \Z} (-1)^m
\frac{q^{\frac{m}{2} (3m+1) +ms  }}{1-q^{m +\frac{a}{c}  }}
+ q^{ \frac{3a}{2c}}
\sum_{m \in \Z} (-1)^m
\frac{q^{\frac{3m}{2} (m+1) +ms  }}{
\left(1-q^{m +\frac{a}{c}}  \right)^2
}\right).
\end{eqnarray*}
Now Lemma \ref{thetaLemma} implies the following decomposition of the
second rank moment; the subsequent lemma describes the corresponding
theta integral component.

\begin{corollary}
We have
\begin{eqnarray*}
\mathcal{R} \left(\frac{a}{c};-\frac{1}{z} \right)
= \frac{(-i z)^{ \frac{3}{2}  }}{48 \sqrt{6}} \, e^{-\frac{\pi i
z}{288}  + \frac{\pi i s a z}{12c}  -
\pi i \frac{a^2z}{8c^2}}
T_2\left(\frac{a}{c};q^{\frac{1}{24}} \right)
+\frac{(-i z)^{\frac32 }}{64 \sqrt{3}}
\int_{ 0}^{\infty}
\frac{\Theta_{a,c}\left( \frac{iu}{16d_c^2}    \right)}{(   -i
(iu+z)   )^{ \frac32 }} du.
\end{eqnarray*}
\end{corollary}

\begin{lemma}
\begin{eqnarray*}
\mathcal{N}_{  \frac{a}{c}} \left(z+1 \right) &=& \mathcal{N}_{
\frac{a}{c} }(z), \\
\mathcal{ N}_{ \frac{a}{c}} \left( -\frac{1}{z} \right) &=&
- \frac{i}{64\sqrt{3} \pi}  \int_{  -\bar \tau }^{  i \infty }
\frac{\Theta_{ a,c}\left(  \frac{\tau}{16d_c^2} \right)}{(   -i (
\tau +z  )     )^{  \frac32 }}
d\tau
+ \frac{i}{64\sqrt{3} \pi}  \int_{ 0}^{\infty } \frac{\Theta_{
a,c}\left(  \frac{it}{16d_c^2} \right)}{(   -i ( z +it  )     )^{
\frac32 }}\, dt.
\end{eqnarray*}
\end{lemma}
By work of Shimura it now follows  that
$\Theta_{a,c}\left(-\frac{1}{16d_c^2z} \right) (-iz)^{ \frac12}$ is a
modular form of weight $\frac12$ on $\Gamma_1\left( 64
d_c^2\right)$.
Moreover observe that
\begin{eqnarray*}
L \left( e^{-3 \pi i \gamma w^2 z} R_2 \left( \frac{a}{c},q; - z
w\right)  \right)
=  -z R_2 \left( \frac{a}{c};q  \right).
\end{eqnarray*}
That  $\mathcal{M}_{\frac{a}{c}}$  is annihilated under
$\Delta_{\frac{3}{2}}$ can be seen as  in \cite{B1}.
Combining the above now easily gives the theorem.

\section{Proof of Theorems \ref{quasitheorem} and
\ref{maintheorem}}  \label{MainSection}
\subsection{Proof of Theorem \ref{quasitheorem}}
\label{quasiSubSection}
First observe that  the function  $R_{k}(q)$ is a linear combination
of the functions $\mathcal{R}_j(q)$ given in Section
\ref{RCSection}.
We  use the rank-crank relation (\ref{eq:rankcrank})   and the
modularity of the occurring functions.
The functions $C_{\alpha}(q^{24}) q^{-1}$ are quasimodular forms and
$\eta(24z)$ is a modular form.
Since $R_2 \left( q^{24}\right)q^{-1}$ is a quasimock theta function,
it follows inductively that the functions $R_{\alpha}\left(q^{24}
\right)q^{-1}$ are also quasimock theta functions.

For the readers convenience we  give more details   in the case
$k=3$.
We obtain from (\ref{eq:rankcrank})
\begin{eqnarray*}
3 \mathcal{R}_4(q)
= - 2 \left(3 \delta_q + 1  \right) C_2(q) +8 C_4(q)
+  3\left( -12 \delta_q+1\right)  \mathcal{R}_2(q).
\end{eqnarray*}
Since
\begin{eqnarray*}
\eta_2(n) &=& \frac12 N_2(n), \qquad \text{and} \\
\eta_4(n) &=&\frac{1}{24} ( N_4(n) -N_2(n) ),
\end{eqnarray*}
this   implies that
\begin{multline}\mylabel{eq:Maassrelation}
36q^{-1} R_3 \left(q^{24} \right) \\
=- \frac{9}{8} q^{-1} C_2 \left( q^{24} \right)
-\frac{1}{8} \delta_q \left(   q^{-1}  C_2\left(q^{24}\right)
\right)
+ 4  q^{-1}C_4\left(q^{24}\right) -\frac{3}{2} q^{-1}
R_2\left(q^{24}\right)
- \frac{3}{2}   \delta_q \left( q^{-1} R_2\left(q^{24}\right)
\right).
\end{multline}

From Section  \ref{RCSection}, we know   that  the functions $q^{-1}
C_j\left(q^{24}\right)$ are  also quasimodular forms.  This further
implies that the function
$\delta_q\left( q^{-1} C_2(q^{24})\right)$ is as well, since
differentiation preserves the space of quasimodular forms.
Moreover
\begin{eqnarray*}
q^{-1} R_2 \left(q^{24} \right) = \mathcal{R}(q)
= \mathcal{M}(z) + \mathcal{N}(z) +\frac{1}{ 24 \eta(24z)} -
\frac{E_2(24z)}{8 \eta(24z)}.
\end{eqnarray*}
Therefore $q^{-1} R_3 \left(q^{24} \right)$ can be written as the sum
of a quasimodular form and  the derivative of the holomorphic part of
a weak Maass form, and is thus a quasimock theta function.

\subsection{Proof of Theorem  \ref{maintheorem}}
Denote by $\mathcal{D}_k(m_1,m_2, \cdots, m_k;n)$ the number of
$k$-marked Durfee symbols arising from partitions of $n$ with $i$th
rank equal to $m_i$. Let
\begin{eqnarray*}
R_k(z_1,\cdots, z_k;q):=
\sum_{ m_1,\cdots,m_k \in \Z }\,  \sum_{n=0}^{\infty}
\mathcal{D}_k(m_1,m_2, \cdots, m_k;n)\, z_1^{m_1}z_2^{m_2}\cdots
z_k^{m_k}q^n.
\end{eqnarray*}
In particular
\begin{equation*}
R_k(w;q)
=
R_k\left(w,w^2, \cdots, w^k;q \right)
.
\end{equation*}

If  $x_i \not=x_j$  and $x_ix_j\not=1$, then Andrews showed that the
generating function for the Durfee symbols is actually a linear
combination of rank functions \cite{An2}:
\begin{eqnarray}\mylabel{eq:smallK}
R_{k}(x_1,x_2, \cdots, x_k;q)
= \sum_{i=1}^k \frac{R(x_i;q)}{\prod_{\substack{ j=1\\j \not =i }
}^{k}(x_i-x_j) \left(1-\frac{1}{x_ix_j} \right)}.
\end {eqnarray}
Standard techniques for dissections of $q$-series then give
\begin{equation}\mylabel{eq:ortho}
\sum_{n=0}^{\infty} NF_k(r,t;n)\, q^{24n-1}
= \frac{1}{t}
\left(
R_{k}\left(q^{24}\right)
+\sum_{j=1}^{t-1}
\zeta_t^{-rj}  q^{-1}\,R_k \left(\zeta_t^j;q^{24} \right)  \right)
.
\end{equation}
Thus, if    $\zeta_t^{lj} \not = \zeta_t^{mj}$ and
$\zeta_t^{j(l+m)}\not=1$ for all
$0<l \not=m \leq k$, which is guaranteed if $k \leq \frac{p_t}{2}$,
then
we obtain by (\ref{eq:smallK})
\begin{eqnarray}\mylabel{eq:rootunity}
R_k\left(\zeta_t^j;q \right)
= \sum_{l=1}^{k}
\frac{R \left( \zeta_t^{lj};q \right)
}{
\prod_{ \substack{ m=1\\m \not=l    }}^k \left( \zeta_t^{jm}
-\zeta_t^{jl} \right)
\left( 1-\zeta_t^{ -j(l+m) } \right)
}.
\end{eqnarray}

Otherwise, for ``large" $k$ we need a modified version of
(\ref{eq:smallK}). Namely  if
$x_i=x_j$ or $x_ix_j=1$, then $R_k(x_1,x_2, \cdots, x_k;q)$ can be
related to $R(x;q)$ via  analytic continuation of (\ref{eq:smallK}).
One can show  that the  new function is a linear combination of
$\left[ \frac{\partial^r R(y;q)}{\partial y^r} \right]_{y=x_i}$.
For example if $k=2$,  $x_1=x_2$, and $x_1x_2 \not=1$, then
\begin{eqnarray*}
R_2(x_1,x_2;q)
= \lim_{x_2 \to x_1}
\left(\frac{R(x_1;q)}{\left(x_1-x_2\right)\left(1-\frac{1}{x_1x_2}
\right)} +
\frac{R(x_2;q)}{\left(x_2-x_1\right)\left(
1-\frac{1}{x_1x_2}\right)}\right)
= \frac{\left[ \frac{\partial}{\partial y} R(y;q
)\right]_{y=x_1}}{\left(1-x_1^{-2}\right)}.
\end{eqnarray*}
In general, $R_k(x_1, \dots, x_k; q)$ will be a linear combination of
various derivatives of $R(y;q)$ for any values assigned to the
$x_i.$  This can be seen by comparing the two sides of
(\ref{eq:smallK}); since the left side has no poles, all of the
singularities on the right side must be removable, and thus
L'Hospital's rule can be applied (at most $2k$ times).

We largely restrict our argument to the case that  $k\leq
\frac{p_t}{2}$, and then make some comments on the general case.
 Fix a prime $p>3$ with $p \nmid t$.
 We treat the two summands in (\ref{eq:rootunity}) separately.
 From Subsection  \ref{quasiSubSection}, we know that
 $R_{\alpha}\left( q^{24}\right)q^{-1}$ is a quasimock theta
 function.  Moreover one can conclude from Theorem \ref{maintheorem2}
 that the holomorphic part of the associated weak Maass form is
 supported on negative squares. Thus the restriction to coefficients
 lying in
 $$
 S_p:= \left\{ n \in \Z: \leg{ 24   n - 1}{p} = - \leg{- 1}{p}
 \right\}
 $$
 is a  weakly holomorphic modular form on $\Gamma_1 \left( 96
 d_t^2p^2\right)$.
 Moreover from the work of the first  author and Ono \cite{BO2}, we
 know that  the restriction of $R \left(\zeta_t^{j};q^{\ell_t  }
 \right)q^{ -\frac{\ell_t}{24}}$ to
 those coefficients lying in  $S_p$
 is a weakly holomorphic  modular form on $\Gamma_1
 \left(6f_t^2l_tp^2 \right)$.
 Finally, work of Serre implies that quasimodular forms are actually
 $p$-adic modular forms, and the proof concludes as in \cite{B1}.

To prove the general case, we first consider the function $\left[
\frac{\partial}{\partial w} R(w;q)\right]_{w= \zeta_c^a }$ which
is again the base case for induction. We have
\begin{equation*}
\left[ \frac{\partial}{\partial w} R(w;q)\right]_{w= \zeta_c^a }
=
\frac{1}{(q;q)_{\infty} }\sum_{ n \in \Z}
\frac{(-1)^{n+1} q^{ \frac{n}{2}(3n+1)   }}{\left( 1-\zeta_c^aq^n
\right)}
- \frac{\zeta_c^a\left(1-\zeta_c^a \right)}{(q;q)_{\infty}}
\sum_{ n \in \Z}
\frac{(-1)^{n+1} q^{ \frac{3n}{2}(n+1)   }}{\left( 1-\zeta_c^aq^n
\right)^2}.
\end{equation*}
Since we know that $R\left(\frac{a}{c};q \right)$ is related to a
weak Maass form, it is enough to consider the function  $R_2
\left(\frac{a}{c};q \right)$ defined in (\ref{eq:TwistGen}).
Theorem \ref{maintheorem2} implies that $\mathcal{R}
\left(\frac{a}{c};q \right)$ is the holomorphic part of a weak Maass
form. Moreover as in \cite{BO2}, one can compute that the
non-holomorphic part is supported on negative squares. Now we argue
as before, utilizing the easy shown technical result that for these
Maass forms, the quadratic twist operator commutes with
differentiation (due to the form of the non-holomorphic coefficients).
For higher derivatives, it is enough to relate the functions
$\mathcal{R}_{j,a,c}(q)$ to quasiweak Maass forms. For this we apply
$\partial_w^j$ to (\ref{eq:rankcrank}), set $w=\zeta_c^a$, and then
argue inductively.
The case $k=2$ again begins the  induction.
The claim now follows, using the relation of  the functions $C_{j,a,c}(q)$
to quasimodular forms, and the modularity of the twisted Eisenstein
series $G_{j+1}^{(0,a)}(z)$.

\section{Explicit differential operators and rank moment differences}
\label{CongruenceSection}

In this section, we show more explicitly how the modularity and holomorphicity of the higher rank moments $\mathcal{R}_{2k}$ are determined by the derivatives of $\mathcal{R}_2.$  We describe the differential operator
that naturally acts on $\mathcal{R}_2$, and then prove Theorem \ref{Congruencetheorem} by considering
the operator $\ell$-adically.
Throughout this section we assume that $\ell>3$ is prime.

\subsection{Higher rank moments and derivatives of $\mathcal{R}_2$}
We begin by defining what we mean by an
``$\ell$-integral quasimodular form of level $1$.''
We consider functions $E_{2}^{a}(z) \, F_b(z)$, where $F_b(z)\in M_b(1)$,
the coefficients in the $q$-expansion
of $F_b(z)$ are $\ell$-integral and have bounded denominators, and
$a$ and $b$ are nonnegative integers. We call such a function
an $\ell$-integral quasi-modular form of weight $2a+b$.
Let $k$ be a nonnegative integer. In general, an
\textit{$\ell$-integral quasi-modular form of weight $k$ and level $1$} is sum of such functions
where $2a+b=k$.
For $k$ an even nonnegative integer, let $\mathcal{X}_k$ denote
the set of functions that are sums of $\ell$-integral quasi-modular forms
of weight $\le k$. Let
\beqs
P\mathcal{X}_k = \{ G\,P\,:\,G\in\mathcal{X}_k\}.
\eeqs

In equation \eqn{rankcrank}, replace $a$ by $2k$ and denote the left
side by $Y_{2k}$, so that
\beq
Y_{2k} :=
\sum_{i=0}^{k-1}\binom{2k}{2i}
\sum_{\substack{\alpha+\beta+\gamma=2k-2i\\
      \alpha, \beta, \gamma \; \text{even} \; \ge 0}}
\binom{2k-2i}{\alpha,\beta,\gamma}\,C_\alpha\,C_\beta\,C_\gamma\,P^{-2}
- 3\left(2^{2k-1}-1\right) C_2.
\mylabel{eq:Ydef}
\eeq
With this notation, the relation \eqn{rankcrank} can be solved for $\mathcal{R}_{2k},$ yielding
\begin{align}
\mathcal{R}_{2k} &= \frac{1}{(2k-1)(k-1)}Y_{2k} -
         \frac{1}{(2k-1)(k-1)}\left(
\sum_{i=1}^{k-1} 6\binom{2k}{2i}\left(2^{2i-1}-1\right)
\delta_q(\mathcal{R}_{2k-2i})\right. \mylabel{eq:Rrec}\\
&\qquad + \left.\sum_{i=1}^{k-1}\left[
\binom{2k}{2i+2}\left(2^{2i+1}-1\right) - 2^{2i}\binom{2k}{2i+1}
+ \binom{2k}{2i}
\right] \mathcal{R}_{2k-2i}\right).
\nonumber
\end{align}
An easy induction argument shows that the differential operator that acts on $\mathcal{R}_2$ is a polynomial $P_k(\delta_q)$ with rational coefficients and degree $k-1,$ so that
\beq
\mathcal{R}_{2k} = P_k(\delta_q)\, \mathcal{R}_2 + \sum_{j=2}^k Q_{k,j}(\delta_q)\,Y_{2j},
\mylabel{eq:Rsolprop}
\eeq
where $Q_{k,j}\in\Q[x]$ has degree $k-j$.

We now focus on $P_k(x)$; we will see shortly that the other terms in (\ref{eq:Rsolprop}) may be absorbed into
the quasimodular component of $\mathcal{R}_{2k}.$
\begin{prop}
\mylabel{prop:Pkrec}
Let $P_0(x) := 0$ and $P_1(x) := 1.$  For $k\ge2$ we have the recurrence relation
\begin{equation}
P_k(x) = (1 - 12x)\,P_{k-1}(x) - 36x^2\,P_{k-2}(x),
\notag
\end{equation}
and the explicit formula
\begin{equation}
P_k(x) = 2^{1-2k}\sum_{j=0}^{k-1}\binom{2k}{2j+1} (1 - 24x)^j\mylabel{eq:Pk1}
       = \frac{1}{\sqrt{1-24x}} \left(
\left(\frac{1+\sqrt{1-24x}}{2}\right)^{2k}
-
\left(\frac{1-\sqrt{1-24x}}{2}\right)^{2k}
\right).
\notag
\end{equation}
\end{prop}
\begin{proof}
From \eqn{Rrec} and \eqn{Rsolprop} we see that for $k \geq 0,$
\begin{align}
P_k(x) &= -\frac{1}{(2k-1)(k-1)}
\sum_{i=1}^{k-1} \left(6x\binom{2k}{2i}\left(2^{2i-1}-1\right)
+ \binom{2k}{2i+2}\left(2^{2i+1}-1\right)\right)\notag \\
&\qquad\qquad\qquad\qquad\qquad\qquad\qquad    \left. - 2^{2i}\binom{2k}{2i+1} + \binom{2k}{2i}\right)
P_{k-i}(x).
\nonumber
\end{align}
For $k\ge0$ we define
\beq
V_k (z) := P_k\left(\frac{1-z^2}{24}\right).
\notag
\eeq
Using some elementary (though lengthy) binomial sum evaluations and an induction
argument, one can show that
\beq
V_k(z) = \frac{1}{z} \left(
\left(\frac{1+z}{2}\right)^{2k}
-
\left(\frac{1-z}{2}\right)^{2k}
\right).
\notag
\eeq
From this one can conclude that for $k \geq 2,$
\beq
V_k(z) = \left(\frac{1+z^2}{2}\right)\,V_{k-1}(z)
         -\left(\frac{1-z^2}{4}\right)^2\,V_{k-2}(z).
\notag
\eeq
Letting $z=\sqrt{1-24x}$ (taking any fixed branch of the square root), we obtain
the results.
\end{proof}

By letting $2k=\ell+1$ and reducing the formula from Proposition \ref{prop:Pkrec} modulo $\ell$ we obtain the
following result.
\begin{corollary}
\mylabel{cor:Pcong}
For $k\ge2$ the polynomial $P_k(x)$ has integer coefficients.
If $\ell>3$ is prime then
\beq
P_{\tfrac{\ell+1}{2}}(x) \equiv \frac{\ell+1}{2}
\left(1 + (1 -24x)^{\tfrac{\ell-1}{2}}\right)
\pmod{\ell}.
\notag
\eeq
\end{corollary}
We will need precise expansions of the rank and crank moments into $\ell$-integral components.
\begin{prop}
\mylabel{propo:diffint}
Let $\ell>3$ be prime.
\begin{enumerate}
\item 
For $k\ge0$ even we have
$\delta_q\left(\mathcal{X}_k\right) \subset \mathcal{X}_{k+2}$.
\item 
For $k\ge0$ even and $m\ge0$ we have
$\delta_q^m\left(P\mathcal{X}_k\right) \subset P\mathcal{X}_{k+2m}$.
\item 
For $1 \le j \le \tfrac{\ell+1}{2}$ except $j=\frac{\ell-1}{2}$
we have $C_{2j} \in P\mathcal{X}_{2j}$.
\item 
$C_{\ell-1} = 2 P\Phi_{\ell-2} + P G$
for some $G\in\mathcal{X}_{\ell-1}$.
\item  
For $1 \le j \le \tfrac{\ell-3}{2}$
we have $Y_{2j} \in P\mathcal{X}_{2j}$.
\item 
$Y_{\ell-1} = 6 P\Phi_{\ell-2} + P G$
for some $G\in\mathcal{X}_{\ell-1}$.
\end{enumerate}
\end{prop}
\begin{proof}
Suppose $\ell>3$ is prime and $k\ge0$ is even.
Our goal is to write the crank moments $C_{2j}$ and hence the $Y_{2j}$
in terms of $\ell$-integral quasimodular forms.
Throughout the proof, we will only
write the components that are {\it not} obviously $\ell$-integral,
with trailing ellipses representing the remaining terms which are
$\ell$-integral quasimodular forms.

(1) \quad This result is well known. It follows from the fact that
$$
12 \delta_q E_2 =E_2^2  - E_4\quad\mbox{and}\quad
12 \delta_q F - k E_2 F \in M_{k+2} \quad \mbox{if $F\in M_k$}.
$$

(2) \quad
From \cite[(2.11)]{AG} we have
\beq
\delta_q P = \Phi_1 P = \tfrac{1}{24}(1 - E_2)P \in P\mathcal{X}_2.
\mylabel{eq:dqP}
\eeq
The result follows from (1) and \eqn{dqP} by an induction argument.

(3) \quad Suppose $1\le j \le \tfrac{\ell-3}{2}$. Then
$\Phi_{2j-1}=\frac{B_{2j}}{4j}(E_{2j}-1)\in \mathcal{X}_{2j}$
by the von-Staudt and Kummer congruences \cite{BS}, \cite[p.20]{SWD}.
Hence
$C_{2j} \in P\mathcal{X}_{2j}$
by \eqn{relation}. Similarly,
$C_{\ell+1} \in P\mathcal{X}_{\ell+1}$
by \eqn{relation}, since
\begin{align*}
C_{\ell+1} &= 2\Phi_{\ell} P + 2 \ell C_{\ell-1} \Phi_1 + \cdots
+  2 \Phi_{\ell -2} C_2, \\ 
           &= 2\Phi_{\ell} P + 8 \ell \Phi_{\ell-2} \Phi_1 + \cdots
\nonumber
\end{align*}
and
\beq
\ell \Phi_{\ell-2}\Phi_1
= \frac{\ell B_{\ell-1}}{2(\ell-1)}(1 - E_{\ell-1})\Phi_1
\in \mathcal{X}_{\ell+2},
\notag
\eeq
again by the von-Staudt and Kummer congruences.

(4) \quad The result follows from \eqn{relation} and (3).

(5) \quad The result follows from \eqn{Ydef} and (3).

(6) \quad The result follows from \eqn{Ydef}, (3), (4)
and since
\beqs
Y_{\ell-1}  = 3 C_{\ell -1} + \cdots = 6 \Phi_{\ell-2} P + \cdots.
\eeqs
\end{proof}

We can now prove that the $\ell$-adic behavior of the higher rank moments comes from that
of $\mathcal{R}_2.$
\begin{theorem}
\mylabel{thm:bigR2}
Let $\ell>3$ be prime.
Then
\beq
\mathcal{R}_{\ell+1} - P_{\tfrac{\ell+1}{2}}(\delta_q) \mathcal{R}_2
\in P \mathcal{X}_{\ell+1}.
\notag
\eeq
\end{theorem}
\begin{proof}
The idea of the proof is to use \eqn{Rsolprop}, and rewrite the
$Y_{2j}$ in terms of $\ell$-integral quasi-modular forms,
keeping track of when $\ell$ occurs in denominator of a coefficient.
We consider the equations \eqn{Rrec} for $2 \le k \le \frac{\ell+1}{2}$.
In these equations the only time $\ell$ occurs in a denominator
is when $k=\frac{\ell+1}{2}$, and only the term of the right side
which is not $\ell$-integral is the term involving $Y_{\ell+1}$.
We note that all coefficients in the sum are $\ell$-integral.
The only term which could rise to a non-$\ell$-integral quasi-modular
form is the term with $i=1$.
Equation \eqn{Ydef} gives $Y_{\ell+1}$ in terms of crank moments.
We wish to write $\tfrac{1}{\ell} Y_{\ell+1}$ in terms of quasi-modular
forms identifying terms that are not $\ell$-integral.
As usual,
we will only write the components that are {\it not} obviously $\ell$-integral,
with trailing ellipses representing the remaining $\ell$-integral portion.
We find that
\begin{align}
\tfrac{1}{\ell} Y_{\ell+1} &= \tfrac{3}{\ell}\left(C_{\ell+1}-C_2\right)
+ 3 (\ell + 1) C_{\ell-1} C_2 P^{-1} + \tfrac{3}{2} (\ell + 1) C_{\ell-1}
+ \cdots
\nonumber\\
& = \tfrac{6}{\ell} P\left(\Phi_{\ell} - \Phi_1\right)
 + 18 (\ell + 1) \Phi_{\ell-2} \Phi_{1} P + 3 (\ell+1) \Phi_{\ell-2} P
 + \cdots.
\nonumber
\end{align}
Therefore by \eqn{Ydef}, \eqn{Rrec}, \eqn{Rsolprop},
and Proposition \propo{diffint}
\begin{align}
\mathcal{R}_{\ell+1} - P_k(\delta_q) \mathcal{R}_2
& =  \frac{2}{\ell(\ell-1)} Y_{\ell+1}
- \frac{ ( \ell+1 )  ( 72\delta_q+7{\ell}^{2}-37\ell+42 ) }
       {6 ( \ell-1 )  ( \ell-2 )  ( \ell-3 )}\, Y_{\ell-1} + \cdots
\mylabel{eq:R2}\\
& = F_1 +
\sum_{j=2}^{(\ell+1)/2} \Qtwid_{\ell,j}(\delta_q) P Z_{2j},
\nonumber
\end{align}
where the $\Qtwid_{\ell,j}$ are polynomials of degree $\ell-j$ with
integer coefficients, the $Z_{2j}\in\mathcal{X}_{2j}$
($2 \le j \le \frac{\ell+1}{2}$),
and
\begin{align}
F_1 &= \frac{2}{\ell-1}\left(
\tfrac{6}{\ell} P\left(\Phi_{\ell} - \Phi_1\right)
 + 18 (\ell + 1) \Phi_{\ell-2} \Phi_{1} P + 3 (\ell+1) \Phi_{\ell-2} P
\right)
\nonumber\\
&\qquad - \frac{ ( \ell+1 )  ( 72\delta_q+7{\ell}^{2}-37\ell+42 ) }
       {6 ( \ell-1 )  ( \ell-2 )  ( \ell-3 )}\,
   \left(
   6 \Phi_{\ell-2} P
   \right)
\nonumber\\
& = F_2 P + \frac{1}{\ell} F_3,
\nonumber
\end{align}
where
$F_2\in\mathcal{X}_{\ell+1}$
and
after some calculation we find that
\beq
F_3=
    (\tfrac{1}{2} - 12 \delta_q) \Etwid_{\ell-1} P
   + 12 \Etwid_{\ell+1} P
    - \tfrac{3}{2} E_2 \Etwid_{\ell-1} P.
\notag
\eeq
Here we have defined the $\ell$-integral modular forms
$\Etwid_{\ell-1}(z):=\frac{\ell B_{\ell-1}}{2(\ell-1)} - \ell
\Phi_{\ell-2}$ and $\Etwid_{\ell+1}(z):=\frac{B_{\ell+1}}{2(\ell+1)}
- \Phi_{\ell}$,  and we have used the
von-Staudt and Kummer congruences $\ell B_{\ell-1} \equiv -1 \pmod{\ell}$,
and $12 B_{\ell+1} \equiv 1 \pmod{\ell}$.
Since $\Etwid_{\ell-1}\in M_{\ell-1}(1)$,
\beq
V_{\ell} = 12 \delta_q \Etwid_{\ell-1} - (\ell -1) E_{2}\Etwid_{\ell-1}
\notag
\eeq
is an $\ell$-integral modular form of weight $(\ell+1)$.
Now
\beqs
12 \delta_q \Etwid_{\ell-1}P  = 12\left(\delta_q \Etwid_{\ell-1}\right) P
+ 12 \Etwid_{\ell-1} \delta_q P
= V_\ell P + (\ell - \tfrac{3}{2}) E_{2}\Etwid_{\ell-1} P
+ \tfrac{1}{2}\Etwid_{\ell-1}P.
\eeqs
Hence
\beq
F_3 = (-V_\ell  + 12 \Etwid_{\ell+1})P - \ell E_2\Etwid_{l-1} P.
\notag
\eeq
Now
\beqs
-V_\ell  + 12 \Etwid_{\ell+1} \equiv  -E_2 \Etwid_{\ell-1}
+ 12\Etwid_{\ell+1} \equiv 0 \pmod{\ell},
\eeqs
by the well known congruences $2\Etwid_{\ell-1}\equiv1$ and
$24\Etwid_{\ell+1}\equiv E_2 \pmod{\ell}$,
and we note that $-V_\ell+12 \Etwid_{\ell+1}
\in M_{\ell+1}(1)$.  Hence
\beq
F_1 = F_4 P,
\notag
\eeq
where $F_4\in\mathcal{X}_{\ell+1}$.
The function
\beqs
\sum_{j=2}^{(\ell+1)/2} \Qtwid_{\ell,j}(\delta_q) P Z_{2j}
\eeqs
has the same property by Proposition \propo{diffint}, and the result
follows from \eqn{R2}.
\end{proof}

For $\ell>3$ prime and  $\epsilon\in\{-1,0,1\}$ we define the operator
$U_{\epsilon,\ell}^{*}$,
which acts on $q$-series by
\beq
U_{\epsilon,\ell}^{*}\left(\sum_n a(n) q^n \right)
:= \sum_{\leg{1-24n}{\ell}=\epsilon} a(n) q^n.
\notag
\eeq
The following corollary of Theorem \thm{bigR2} follows
from Corollary \ref{cor:Pcong}.
\begin{corollary}
\mylabel{cor:bigR2cor}
Let $\ell>3$ be prime and suppose $\epsilon= -1$ or $0$.
Then
\beq
U_{\epsilon,\ell}^{*}\left( \mathcal{R}_2 \right)
\equiv
U_{\epsilon,\ell}^{*}\left( G_{\ell} P\right) \pmod{\ell},
\notag
\eeq
where
$G_{\ell}\in\mathcal{X}_{\ell+1}$.
\end{corollary}
\begin{proof}
Let $\epsilon=-1$ or $0$.
We define
\beq
c_{\epsilon,\ell} :=
\begin{cases}
\frac{\ell+1}{2} & \mbox{if $\epsilon=0$},\\
1 & \mbox{if $\epsilon=-1$}.
\end{cases}
\notag
\eeq
If $\leg{1-24n}{\ell}=\epsilon$
we have
\beq
\left(1- P_{\tfrac{\ell+1}{2}}(n)\right) \equiv c_{\epsilon,\ell}
\pmod{\ell},
\mylabel{eq:Pkmodep}
\eeq
by Corollary \corol{Pcong}.
Since $\mathcal{R}_{\ell+1}\equiv\mathcal{R}_2\pmod{\ell}$,
from Theorem \thm{bigR2} we have
\beq
U_{\epsilon,\ell}^{*}\left(\left(1- P_{\tfrac{\ell+1}{2}}(\delta_q)\right) \mathcal{R}_2
\right)
\equiv  U_{\epsilon,\ell}^{*}\left(G_{\ell} P\right) \pmod{\ell},
\notag
\eeq
where
$G_{\ell}\in\mathcal{X}_{\ell+1}$.
Finally by \eqn{Pkmodep} we have
\begin{align}
U_{\epsilon,\ell}^{*}\left( \left(1- P_{\tfrac{\ell+1}{2}}(\delta_q)\right) \mathcal{R}_2
\right)
& =
\sum_{\leg{1-24n}{\ell}=\epsilon} \left(1- P_{\tfrac{\ell+1}{2}}(n)\right) N_2(n)
q^n
\nonumber\\
& \equiv  c_{\epsilon,\ell}
\sum_{\leg{1-24n}{\ell}=\epsilon} N_2(n) q^n
\equiv c_{\epsilon,\ell} U_{\epsilon,\ell}^{*}\left( \mathcal{R}_2 \right) \pmod{\ell},
\nonumber
\end{align}
and the result follows since $c_{\epsilon,\ell}\not\equiv0\pmod{\ell}$.
\end{proof}

The following theorem describes the relation between $\mathcal{R}_{2k}$ and $\mathcal{R}_2$ for other $k$.
The proof is entirely analogous to Theorem \thm{bigR2}, with an extra term that is not $\ell$-integral
appearing when $k = (\ell-1)/2.$
\begin{theorem}
\mylabel{thm:RkCkbig}
Let $\ell > 3$ be prime.
\begin{enumerate}
\item
\label{thm:RkCkbigC}
For $1 \le k \le \tfrac{\ell-3}{2}$ we have $C_{2k} \in P \mathcal{X}_{2k},$ whereas $C_{\ell-1}
= 2P \Phi_{\ell-2} + P G_C$ for some $G_C\in\mathcal{X}_{\ell-1}$.

\item
\label{thm:RkCkbigR}
For $2 \le k \le \tfrac{\ell-3}{2}$ we have
\begin{equation}
\mathcal{R}_{2k} - P_k(\delta_q) \mathcal{R}_2 \in P \mathcal{X}_{2k},
\notag
\end{equation}
whereas
\begin{equation}
\mathcal{R}_{\ell-1} - P_{\tfrac{\ell-1}{2}}(\delta_q) \mathcal{R}_2
= 2P \Phi_{\ell-2} + P G_R
\notag
\end{equation}
for some $G_R \in \mathcal{X}_{\ell-1}.$

\end{enumerate}
\end{theorem}

\subsection{Proof of Theorem \ref{Congruencetheorem}}
We are now ready to prove Theorem \ref{Congruencetheorem} by combining the above results
with the theory of $\ell$-adic modular forms.

\noindent (1) \quad
By Corollary \corol{bigR2cor},
\beq
U_{0,\ell}^{*}(\mathcal{R}_2) \equiv U_{0,\ell}^{*}(G_\ell P)
\pmod{\ell}
\notag
\eeq
for some $G_\ell\in\mathcal{X}_{\ell+1}$. By reduction modulo ${\ell}$
we may assume that all forms involved have integer coefficients.
The function $G_\ell P$ is a sum of functions of the form
\beq
\Etwid_2^{a} F P \equiv \Etwid_{\ell+1}^a F P
\equiv q^{\tfrac{1}{24}}\frac{\Ftwid_{a}(z)}{\eta(z)}\pmod{\ell},
\notag
\eeq
where $\Etwid_2=\tfrac{1}{24}E_2$ and
$F(z)$ is an integral modular form $F\in M_{b}(1)$,
with total weight $2a + b \le \ell + 1$. Define the associated
$\ell$-adic modular form
$\Ftwid_{a,b}(z):=\Etwid_{\ell+1}^a(z) F(z) \in M_{b+a(\ell+1)}(1)$,
which has weight at most $\tfrac12(\ell+1)^2$.
Define the coefficients $p\left(\Ftwid, n\right)$
by $\sum_n p\left(\Ftwid,\ell n + \beta_\ell\right) q^{24n+r_\ell} =
\frac{\Ftwid(q)}{(q;q)_\infty}$.
By standard arguments using Hecke operators (see \cite{Gar}), we have
\begin{equation}
\sum_{n=0}^\infty N_2(\ell n + \beta_\ell) q^{24n+r_\ell}
\equiv
\sum_{2a + b \le \ell+1}
\sum_{n=0}^\infty p\left(\Ftwid_{a,b},\ell n + \beta_\ell\right) q^{24n+r_\ell}
\equiv \eta^{r_\ell}(24z) G_{\ell,2}(24z) \pmod{\ell},
\notag
\end{equation}
where $G_{\ell, 2} \in \mathcal{X}_{(\ell^2 + 3\ell - r_\ell - 1)/2}.$

\noindent (2) \quad  This claim follows immediately using the same arguments as above and
the relation between $\mathcal{R}_{2k}$ and $\mathcal{R}_2$ found in Theorem \ref{thm:RkCkbig}.

\noindent (3) \quad
To find $\ell$-adic modular forms associated with $\mathcal{R}_{\ell-1}$,
we again use Theorem \ref{thm:RkCkbig}.
The difference is that there is now an extra term that is not
$\ell$-integral, namely
$2 P \Phi_{\ell-2}
= \frac{1}{\ell}
\left(\frac{\ell B_{\ell-1}}{(\ell-1)} - 2\Etwid_{\ell-1} \right) P$.

Define $\alpha_\ell$ such that $\alpha_\ell \equiv \frac{\ell B_{\ell-1}}{(\ell-1)} \pmod{\ell^2},$
so that the first part of the contribution from $2P \Phi_{\ell-2}$ is
\beq
\alpha_\ell \sum_{n=0}^\infty p(n) q^{24n-1} \equiv \alpha_{\ell} \frac{1}{\eta(24z)}
\equiv \alpha_{\ell} \frac{1}{\eta(24z)}
\left( \frac{\eta^\ell(24z)}
            {\eta(24\ell z)} \right)^\ell
\equiv \alpha_\ell \frac{\Delta^{\frac{\ell^2-1}{24}}(24z)}
                        {\eta^\ell(24\ell z)}
\pmod{\ell^2}.
\notag
\eeq
This implies that
\beq
\alpha_\ell\sum_{n=0}^\infty p(\ell n+\beta_{\ell}) q^{24n+r_\ell}
\equiv \alpha_\ell \frac{\Delta^{(\ell^2-1)/24}(24z)\mid U(\ell)}
                        {\eta^\ell(24z)}
\pmod{\ell^2}.
\notag
\eeq
Observe that $\Delta^{(\ell^2-1)/24}(z)\in S_{\tfrac{1}{2}(\ell^2-1)}(1)$, and
recall that $T(\ell) \equiv U(\ell) \pmod{\ell},$ so
\beq
f_\ell(z) := \alpha_\ell \Delta^{(\ell^2-1)/24}(z)\mid T(\ell)
= (q^{(\ell^2-1)/24} + \cdots) \mid T(\ell) = c_3 q^{\lambda_\ell} +
\cdots,
\notag
\eeq
where $\lambda_\ell := \frac{\ell^2 + 24\beta_\ell - 1}{24\ell}$, and $c_3$ is an integer.
Since Hecke operators preserve spaces of modular forms, it must be that $f_\ell(z) = \Delta^{\lambda_\ell}(z) H_1(z)$
for some $H_1(z)\in M_{\tfrac{1}{2}(\ell(\ell-1)-r_\ell-1)}(1)$.  We conclude that
\beq
\alpha_\ell\sum_{n=0}^\infty p(\ell n+\beta_{\ell}) q^{24n+r_\ell}
\equiv \frac{ \Delta^{\lambda_\ell}(24z) H_1(24z)}
            {\eta^\ell(24z)}
\equiv \eta^{r_\ell}(24z) H_1(24z) \pmod{\ell^2}.
\notag
\eeq

We proceed in a similar fashion for the term $2\Etwid_{\ell-1}P.$
Define the coefficients $e_{\ell-1}(n)$ so that
$2\Etwid_{\ell-1,\ell}P \equiv \sum_{n=0}^\infty e_{\ell-1}(n) q^n
\pmod{\ell^2}.$
As before, we find that
\beq
\sum_{n=0}^\infty e_{\ell-1}(n) q^{24n-1}
\equiv 2\frac{\Delta^{(\ell^2-1)/24}(24z) \Etwid_{\ell-1}(24z)}
                        {\eta^\ell(24\ell z)}
\pmod{\ell^2},
\notag
\eeq
and
\begin{align}
\sum_{n=0}^\infty e_{\ell-1}(\ell n+\beta_{\ell}) q^{24n+r_\ell}
&\equiv 2
\frac{\Delta^{(\ell^2-1)/24}(24z) \Etwid_{\ell-1}(24z)\mid U(\ell)}
     {\eta^\ell(24z)}
\pmod{\ell^2}
\notag\\
&\equiv \eta^{r_\ell}(24z) H_2(24z) \pmod{\ell^2},
\nonumber
\end{align}
where
$H_2(z)\in M_{\tfrac{1}{2}(\ell(\ell+1)-r_\ell-3)}(1)$
and
\begin{equation}
\Delta^{\lambda_\ell}(z) H_2(z)
\equiv 2\Delta^{(\ell^2-1)/24}(z) \Etwid_{\ell-1}(z) \mid U(\ell)
\pmod{\ell^{2}}.
\nonumber
\end{equation}
Since $2\Etwid_{\ell-1} \equiv \alpha_\ell \equiv 1 \pmod{\ell},$ we have $H_1(24z) \equiv H_2(24z) \pmod {\ell},$ so the overall
the contribution to the congruence for $U_{0,\ell}^{*}(\mathcal{R}_{\ell-1})$ is
\beq
\frac{1}{\ell} \eta^{r_\ell}(24z)
\left(
H_1(24z) - H_2(24z)
 \right)
\notag
\eeq
as claimed, completing the proof.

\section{Proof of Theorem \ref{IdentityTheorem}} \label{IdentSection}
Using (\ref{eq:ortho}), we observe that
\begin{eqnarray*}
\sum_{n=0}^{\infty}
\left( NF_k(r,t;n) -NF_k(s,t;n)\right)q^{ 24n-1}
=\frac{1}{t}
\sum_{ j=1}^{t-1}
\left( \zeta_t^{ -rj} - \zeta_t^{-sj} \right) q^{-1}R_k
\left(\zeta_t^j;q^{ 24} \right).
\end{eqnarray*}
Without loss of generality, we assume that  $k\leq\frac{p_t}{2}$, the
general case is proven similarly.
In this case we may use (\ref{eq:rootunity}).
The functions $R \left(\zeta_t^j;q^{\ell_t} \right)q^{
-\frac{\ell_t}{24} }$ are the holomorphic parts of weak Maass forms
on $\Gamma_1\left(576t^4\right)$.
One can generalize the usual Atkin $U(t^2)$-operator to weak Maass
forms. This gives   that $q^{-1}R \left(\zeta_t^j;q^{24} \right)$
are the holomorphic parts of weak Maass forms on
$\Gamma_1\left(576t^4\right)$. Moreover by Theorem
\ref{MaassformTheorem}, the
non-holomorphic parts of those forms are supported on  negative
squares.
Generalizing the theory of twists of modular forms to twists of weak
Maass forms, we obtain that  the restriction of those form to the
coefficients  supported on arithmetic progression congruent to $d$
modulo $t$ satisfying  $\leg{1-24d}{t} =-1$   is  a weakly
holomorphic modular form on $\Gamma_1 \left( 576t^6\right)$.
Thus (\ref{Thm:R-1}) follows.

In order to conclude (\ref{Thm:R1}), we have to show that  the restriction of the function
\begin{eqnarray*}
\sum_{j=1}^{t-1}
\left(
\zeta_t^{-rj} - \zeta_t^{-sj}
\right)
\frac{\zeta_t^{2j}}{(1-\zeta_t^j)(\zeta_t^{ 3j} -1)}
\left(q^{-1} R \left(\zeta_t^{j};q^{24} \right)  - q^{-1}R
\left(\zeta_t^{2j};q^{24} \right)
\right)
\end{eqnarray*}
to those arithmetic progressions stated in the theorem
doesn't have a non-holomorphic part. The correct group follows as in
(\ref{Thm:R-1}).
Using Theorem \ref{MaassformTheorem}, we see that this is equivalent
to the identity
\begin{equation}\mylabel{eq:rootidentity}
0= \sum_{ j=1}^{t-1}
\left(\zeta_t^{-rj} -\zeta_t^{-sj} \right)
\frac{\zeta_t^{2j}}{(1-\zeta_t^j)(\zeta_t^{ 3j} -1)}
\left(
\sin \left( \frac{\pi j}{t}\right)\sin \left(\frac{3\pi jd}{t}
\right)
- \sin \left( \frac{2\pi j}{t}\right)\sin \left(\frac{6\pi jd}{t}
\right)
\right).
\end{equation}
Since    $\sin (x) = \frac{1}{2i} \left( e^{ix}+e^{-ix} \right)$,
identity (\ref{eq:rootidentity}) is equivalent to
\begin{equation}\mylabel{eq:rootidentity2}\
0= \sum_{ j=1}^{t-1}\frac{\left(\zeta_t^{-rj} -\zeta_t^{-sj}
\right)\left(1-\zeta_t^{3dj}\right)}{ \left(1-\zeta_t^{ 3j}\right)}
\zeta_{2t}^{3j(1-d) }
\left(  1- \left(\zeta_{2t}^j + \zeta_{2t}^{-j} \right)
\left(\zeta_{2t}^{3dj }+ \zeta_{2t}^{-3dj} \right)
\right).
\end{equation}
This is further equivalent to
\begin{multline}\mylabel{eq:rootidentity3}
0= \sum_{ j=1}^{t-1}\left(\zeta_t^{-rj} -\zeta_t^{-sj} \right)
\zeta_{2t}^{3j(1-d) } \\
\left(1+ \zeta_t^{3j} + \cdots +\zeta_t^{3j(d-1) } \right)
\left(1- \zeta_{2t}^{ j(3d+1)}  - \zeta_{2t}^{ j(-3d+1)} -
\zeta_{2t}^{ j(3d-1)} - \zeta_{2t}^{- j(3d+1)}  \right).
\end{multline}
Identity (\ref{eq:rootidentity3}) can   be verified using the conditions
in the theorem and the standard orthogonality of roots of unity.

 \end{document}